\documentstyle[12pt]{article}

\def\triangulo{\hbox{.\kern2pt.\kern-5pt\raise 4pt\hbox{.}}\kern 5pt}

\def\bbb#1{\hbox {{\gordas #1}}}

\font\gordas = msbm10 at 12pt
\def\bbb#1{\hbox {{\gordas #1}}}
\def\a{{\bbb A}}
\def\erre{{\bbb R}}
\def\ce{{\bbb C}}
\def\o{{\bbb O}}
\def\ze{{\bbb Z}}
\def\ache{{\bbb H}}

\begin{document}
\begin{center}
{\bf The Exponential map on the Cayley-Dickson algebras}\\[1cm]
Guillermo Moreno\\
Departamento de Matematicas. \\
CINVESTAV del I.P.N.\\
Mexico D.F. MEXICO\\
gmoreno@math.cinvestav.mx

\end{center}
\vglue3cm
\noindent
{\bf Abstract:} We study the Exponential map for $\a_n=\erre^{2^n}$, the Cayley-Dickson algebras for $n\geq 1$, which generalize the complex exponential map to Quaternions, Octonions and so forth. As an application, we show that the 
selfmap of the unit sphere in $\a_n,S(\a_n)=S^{2^n-1}$, given by taking $k$-powers has topological degree $k$ for $k$ an integer number, from this we derive a suitable ``Fundamental Theorem of Algebra for $\a_n$.''

\vglue7cm
\noindent
\newpage
\noindent
{\bf Introduction.} The Cayley-Dickson algebras $\a_n=\erre^{2^n}$, for $n\geq 0$ are given by the doubling process of Dickson [1].

For $a,b,x$ and $y$ in $\a_n$ and $\a_{n+1}=\a_n\times \a_n$.
$$(a,b)(x,y)=(ax-\overline{y}b,ya+b\overline{x})$$
and for $x_1,x_2$ in $\a_{n-1},\quad \overline{x}=(\overline{x}_1,-x_2)$ when $x=(x_1,x_2)$.

Thus, if $\overline{x}=x$ in $\a_0=\erre$ then $\a_1=\ce$ the Complex numbers, $\a_2=\ache$ the Quaternionic numbers; $\a_3=\o$ the Octonions numbers and so forth.

As is well known (see [7] and [8]) $\a_n$ is commutative only for $n=0$ and
 $n=1$;   $\a_n$ is associative only for $n=0,1,2$;   $\a_n$ is $alternative$ (i.e., $x^2y=x(xy)$ and $xy^2=x(xy))$ for all $x$ and $y$ in $\a_n$) only if $n=0,1,2,3$.

Also $\a_n$ is $flexible$  (i.e., $x(yx)=(xy)x$ for all $x$ and $y$ in $\a_n$) and $power-associative$  (i.e. $x^mx^k=x^{m+k}$ for all $x$ and $y$ in $\a_n$ and $m$ and $k$ natural numbers) for all $n\geq 0$.

By the classical theorem of Hurwitz: $\a_n$ is $normed$  (i.e., $||xy||=||x||||y||$ for all $x$ and $y$ in $\a_n$) if and only if $n=0,1,2,3$ where $||-||$ denotes the standard norm in $\erre^{2^n}$.

Now $\{e_0=1,e_1,e_2,\ldots,e_{2^n-1}\}$ denotes the canonical basis in $\erre^{2^n}$

$$\a_n=\erre e_0\oplus {\rm Im}(\a_n)$$
where $\erre e_0={\rm Span}\{e_0\}$ are the \underline{real elements} in $\a_n$ and
$${\rm Im} (\a_n):={\rm Span}\{e_1,e_2,\ldots,e_{2^n-1}\}$$ 
are the \underline{pure imaginary} elements in $\a_n$. Thus,for every $x$ in $\a_n$ we have a canonical splitting of $x$ into real
 and imaginary parts: $x=re_0+a$ where $r\in\erre$ and $a\in {\rm Im}(\a_n)$.

Now for all $x$ in $\a_n\;||x||^2=x\overline{x}=\overline{x}x$ and 
$2\langle x,y\rangle=x\overline{y}+y\overline{x}$ where 
$\langle -,-\rangle$ denotes the standard inner product in $\erre^{2^n}$. Thus for $a\in {\rm Im}(\a_n),\quad\overline{a}=-a$ and $a^2=-||a||^2$ and
$$\langle xy,z\rangle=\langle y,\overline{x}z\rangle\quad{\rm and} 
\langle x,yz\rangle=\langle x\overline{z},y\rangle$$ 
for all $x$, $y$ and $z$ in $\a_n$ (see [5]).

The main subject of this paper is the Exponential map
$${\rm exp}(x)=e_0+x+\frac{x^2}{2!}+\ldots\frac{x^m}{m!}+\ldots=\sum^\infty_{m=0}\frac{x^m}{m!}\quad{\rm for}\quad x\in \a_n.$$
In $\S$ 1, we prove that ${\rm exp}(x)$ is well defined, that is, converges for all $n\geq 0$ and that some properties of the complex exponential map are also valid for $n>1$.

In $\S$2 we show that the exponential map is surjective, and that restricted to Im$(\a_n)$ it maps onto $S(\a_n)=S^{2^n-1},$ the unit sphere in $\a_n$; also we prove that the $k$-power map is well defined for $k\in \ze$ and has topological degree $k$.

In $\S$ 3 we use the results of the previous sections to prove the Fundamental Theorem of Algebra (F.T.of A.) for $\a_n\quad n\geq 2,$ which generalizes the F.T. of A. for Quaternions of Eilenberg-Niven [2] that goes back to 1949.

First, we present an $Algebraic$ generalization of the F.T. of A. which only considers polynomials where the variable and the coefficients have linearly dependent imaginary part. This version of the F.T. of A. is a straigthforward generalization of the classical F.T. of A. for $\ce$.

Secondly, we present a $Topological$ generalization of the F.T. of A. for ``polynomials of degree $k$'' which are continuous functions inside the homotopy class of the $k$-power map on $\a_n\cup\{\infty\}=S^{2^n}$ for each $k>0$.

We will show, as well, that the topological F.T. of A. generalizes the Algebraic F. T. of A.
\vglue.5cm
\noindent
{\bf I. Basic Properties.} 
\vglue.5cm
Throughtout this paper, we use extensively that $\erre e_0$ is the Center of $\a_n$ for all $n\geq 1$ and that for $x,y$ and $z$ in $\a_n\quad(xy)z=x(yz)$ if at least one of them is a real element.

\vglue.5cm
\noindent
{\bf Definition:} For $x\in\a_n,{\rm exp}(x)=e_0+x+\frac{x^2}{2!}+\cdots =\sum\limits^\infty_{m=0}\frac{x^m}{m!}$

 where $e_0$ is the unit element in $\a_n$.

Clearly for $x=re_0$ for $r\in \erre$, ${\rm exp}(x)=e^re_0$ where $e^r$ is the usual real
exponent map on $\a_0$. In particular if $0$ in $\a_n$ is the null element
 $${\rm exp}(0)=e^0e_0=1\cdot e_0=e_0$$.
\vglue.5cm
\noindent
{\bf Lemma 1.1} If $a$ is a non-zero element in Im$(\a_n)$ then 
$${\rm exp}(a)={\rm cos}(||a||)e_0+{\rm sin}(||a||)\frac{a}{||a||}.$$
{\bf Proof}: Since $a\in {\rm Im}(\a_n)$ then  $a^2=-||a||^2$ so $a^{2k}=(-1)^k||a||^{2k}$ and 

$a^{2k+1}=a^{2k}a=(-1)^k||a||^{2k}a$ for $k\geq 0$. Therefore 
\begin{eqnarray*}
{\rm exp}(a)&=&\sum^\infty_{m=0}\frac{a^m}{m!}=\sum^\infty_{k=0}\frac{(-1)^k||a||^{2k}}{(2k)!}e_0+\sum^\infty_{k=0}\frac{(-1)||a||^{2k}}{(2k+1)!}a\\
&=&{\rm cos}(||a||)e_0+\frac{1}{||a||}\sum^\infty_{k=0}\frac{(-1)^k||a||^{2k+1}}{(2k+1)!}a\\
&=&{\rm cos}(||a||)e_0+{\rm sin}(||a||)\frac{a}{||a||}
\end{eqnarray*}

\hfill Q.E.D.
\vglue.5cm
\noindent
{\bf Corollary 1.2} For non-zero $s$ in $\erre$ and $a$ non-zero in Im$(\a_n)$ 

$$||{\rm exp}(a)||=1\quad{\rm and}\quad {\rm exp}(sa)={\rm cos}(s||a||)+{\rm sin}(s||a||)\frac{a}{||a||}.$$
{\bf Proof:} $||{\rm exp}(a)||^2={\rm cos}^2(||a||)+{\rm sin}^2(a)\frac{||a||}{||a||}=1$.

Since $\frac{s}{||s||}$ is equal to 1 for $s>0$ and equal to $-1$ for $s<0$ and 
 cosine and sine are even and odd function respectively then
\begin{eqnarray*}
{\rm exp}(sa)&=&{\rm cos}(||sa||)e_0+{\rm sin}(||sa||)\frac{sa}{||sa||}\\
&=&{\rm cos}(|s|||a||)e_0+{\rm sin}(|s|||a||)\frac{sa}{|s|||a||}\\
&=&{\rm cos}(s||a||)e_0+{\rm sin}(s||a||)\frac{a}{||a||}.
\end{eqnarray*}

\hfill Q.E.D.
\vglue.5cm
\noindent
{\bf Theorem 1.3} For $x$ in $\a_n$ and $n\geq 1$, the series exp$(x)$ converges.

If $x=re_0+a$ for $r\in\erre$ and $a$ in Im$(\a_n)$ then

$${\rm exp}(x)=e^r({\rm exp}(a))$$
where $e^r$ is the real exponent map.
\vglue.5cm
\noindent
{\bf Proof:} A direct calculation shows that
\begin{eqnarray*}
e^r({\rm exp}(a))&=&(\sum^\infty_{m=0}\frac{r^m}{m!})(\sum^\infty_{m=0}
\frac{a^m}{m!})\\
&=&(e_0+re_0+\frac{r^2e_0}{2!}+\cdots)(e_0+a+\frac{a^2}{2!}+\cdots)\\
&=&e_0+re_0+a+\frac{r^2e_0}{2!}+ra+\frac{a^2}{2!}+\cdots\\
&=&e_0+(re_0+a)+\frac{(re_0+a)^2}{2!}+\cdots\\
&=&{\rm exp}(x)
\end{eqnarray*}

\hfill Q.E.D.
\vglue.5cm
\noindent
{\bf Corollary 1.4.}For $x$ in $\a_n$, $||{\rm exp}(x)||=e^r$, where $r$ is the real part of $x$.
\vglue.5cm
\noindent
{\bf Proof:} 
$$||{\rm exp}(x)||=||e^r{\rm exp}(a)||=|e^r|||{\rm exp}(a)||=|e^r|=e^r.$$
by Theorem 1.3 and Corollary 1.2.

\hfill Q.E.D.
\vglue.5cm
\noindent
{\bf Example:} The known complex identities $e^{\frac{i\pi}{2}}=i$ and $e^{i\pi}=-1$ correspond in $\a_n$ for $n\geq 1$ to  ${\rm exp}(\frac{a\pi}{2})=a$ and ${\rm exp}(\pi a)= -e_0$ respectively for every $a$ in Im$(\a_n)$ such that $||a||=1$.

Now we show that,the identity  $e^{z+w}=e^z\cdot e^w$ in $\ce$ can be generalized to $\a_n$ for $n\geq 1$, to certain extent.
\vglue.5cm
\noindent
{\bf Definition:} Two elements in $\a_n$ for $n\geq 1$ are \underline{Complex dependent } or \underline{$\ce$-dependent} if their respective pure imaginary parts are linearly dependent.
\vglue.2cm
Notice that, every element in $\a_n$ is $\ce$-dependent with any real element and that, for
 $n=1$ every two elements are $\ce$-dependent, because Im$(\a_1)=\erre e_1$.

Also notice that for every $x$ in $\a_n$, exp$(x)$ and $x$ are $\ce$-dependent.
\vglue.5cm
\noindent
{\bf Lemma 1.5} Let $x$ and $y$ be in $\a_n$ for $n\geq 1$.

(i) If $x$ and $y$ are $\ce$-dependent then $xy=yx$.

(ii) For $n=2$ and $n=3,$ we have that, $x$ and $y$ are $\ce$-dependent if and only if $xy=yx$.
\vglue.5cm
\noindent
{\bf Proof:} First of all, we observe that two elements in $\a_n$ commute if and only if their respective imaginary parts commute.

Let $x=re_0+a$ and $y=se_0+b$ in $\erre e_0\oplus{\rm{Im}}(\a_n)=\a_n.$ So
\begin{eqnarray*}
xy&=&(re_0+a)(se_0+b)=(rse_0+rb+sa+ab)\\
yx&=&(se_0+b)(re_0+a)=(sre_0+sa+rb+ba)
\end{eqnarray*}
Then $xy=yx$ if and only if $ab=ba$.

Suppose that $b=ta$ with $t$ in $\erre$ then
$ab=a(ta)=(ta)a=ba$ and we are done with (1). 

Now notice that if $ab=ba$ for $a$ and $b$ in Im$(\a_n)$ then $(ab)$ is real, because  
$\overline{ab}=\overline{b}\overline{a}=(-b)(-a)=ba=ab$. Moreover, since $2\langle a,b\rangle=a\overline{b}+b\overline{a}=-ab-ba=-2(ab)$ then $\langle a,b\rangle e_0=-ab$.

To prove (ii) recall that for $n=2$ and $n=3$ we have that $a(ab)=a^2b$ so $a(ab)=-a(\langle a,b\rangle e_0)$ implies $a^2b=-||a||^2b=-\langle a,b\rangle a$ and $a$ and $b$ are linearly dependent.

\hfill Q.E.D.
\vglue.5cm
\noindent
{\bf Remarks:} 1.-Notice that the presence of \underline{zero divisors} in $\a_n$ for $n\geq 4$, makes possible to have non-zero elements $a$ and $b$ in Im$(\a_n)$ with $ab=ba=0$ and $a$ orthogonal to $b$ (see [5]).

2.- A further study shows that for $a$ in Im$(\a_n)$ and $n\geq 4,$ the \underline{Centralizer} 
of $a,$ defined by $C_a:=\{b\in {\rm{Im}}(\a_n)|ab=ba\}$ is given by 
$$C_a=\erre a\oplus Ker L_a.$$
where \underline{$Ker L_a$} is the right annihilator of $a$.

3.- A characterization of $\ce$-dependence is given by:For $x$ and $y$ in $\a_n$.

    $x$ and  $y$ are $\ce$-dependent if and only if $x(zy)=(xz)y$ for all 
$z$ in $\a_n$.\\
(See [3], [4] and [6]). (We will not use this).
\vskip.5cm
\noindent
{\bf Lemma 1.6} Let $a$ and $b$ be in Im$(\a_n)$ for $n\geq 2$. If $a$ and 
$b$ are linearly \\
dependent then
$${\rm exp}(a+b)={\rm exp}(a){\rm exp}(b).$$
{\bf Proof:} Notice that if either $a$ or $b$ is null then the assertion is trivial.

Suppose that neither $a$ nor $b$ are null then there is non-zero real number $s$ such that $b=sa$ and 
\begin{eqnarray*}
{\rm exp}(b)&=&{\rm cos}(||sa||)e_0+{\rm sin}(||sa||)\frac{sa}{||sa||}\quad    ({\rm Lemma 1.1}).\\
&=&{\rm cos}(s||a||)e_0+{\rm sin}(s||a||)\frac{a}{||a||}               ({\rm Corollary 1.2})
\end{eqnarray*}
Now by the standard trigonometric identities for addition of angles for cosine and sine we have that
\begin{eqnarray*}
{\rm exp}(a){\rm exp}(b)&=&({\rm cos}(||a||)e_0+{\rm sin}(||a||)\frac{a}{||a||})({\rm cos}(s||a||)+{\rm sin}(s||a||)\frac{a}{||a||})\\
&=&[{\rm cos}(||a||){\rm cos}(s||a||)-{\rm sin}(||a||){\rm sin}(s||a||)]e_0\\
&+&[{\rm cos}(||a||){\rm sin}(s||a||)+{\rm sin}(||a||){\rm cos}(s||a||)]\frac{a}{||a||}\\
&=&[{\rm cos}(||a||+s||a||))e_0+
{\rm sin}(||a||+s||a||))\frac{a}{||a||}\\
&=&({\rm cos}(||a+sa||))e_0+({\rm sin}(||a+sa||))\frac{a+sa}{||a+sa||}
\\
&=&{\rm exp} (a+b)
\end{eqnarray*}
because of $a^2=-||a||^2e_0$ and Corollary 1.2.

\hfill Q.E.D.

\vskip.3cm
\noindent
{\bf Theorem 1.7} If $x$ and $y$ are $\ce$-dependent in $\a_n$ and $n\geq 1$ then
$${\rm exp}(x+y)={\rm exp}(x){\rm exp}(y)$$
{\bf Proof:} Suppose that $x=re_0+a$ and $y=se_0+b$ in $\erre e_0\oplus {\rm Im}(\a_n)$ so $x+y=(r+s)e_0+(a+b)$.
\begin{eqnarray*}
{\rm exp}(x+y)&=&e^{r+s}({\rm exp}(a+b))=e^{r+s}({\rm exp}(a){\rm exp}(b))\\
&=&(e^r{\rm exp}(a))(e^s{\rm exp}(b))\\
&=&{\rm exp}(x){\rm exp}(y).\\
&&\qquad\qquad\qquad
\qquad\qquad\qquad\qquad\qquad\qquad\qquad\qquad\quad{\rm Q.E.D.}
\end{eqnarray*}

\vskip.5cm

\noindent
{\bf Remark:} Theorem 1.7 is the best possible in the following sense:
\vglue.2cm
For $n=1.$ Theorem 1.7 correspond to the exponential law.
\vglue.2cm
For $n\geq 2$. Assume that ${\rm exp}(x+y)={\rm exp}(x){\rm exp}(y)$ for some $x$ and $y$ in $\a_n$ then we are forced to have, at least, the following two conditions: $xy=yx$ and $x(xy)=x^2y$.

Now if $a$ and $b$ are the pure parts of $x$ and $y$ respectively, then $xy=yx$ if and only if $ab=ba$ and $x^2y=x(xy)$ if and only if $a^2b=a(ab)$.(See [6])

If $ab=ba$ then $ab=-\langle a,b\rangle e_0$ so $-||a||^2b=a^2b=a(ab)=-\langle a,b\rangle e_0$ and $a$ and $b$ are linearly dependent. (Notice that $ab\neq 0$ because $a(ab)=a^2b\neq 0$).
\vskip.3cm
\noindent
{\bf Proposition 1.8}. For $x$ in $\a_n$ and $n\geq 1$.

(i) ${\rm exp}(\overline{x})=\overline{{\rm exp}(x)}$.

(ii) ${\rm exp}(-x)={\rm exp}(x)^{-1}$

\vskip.3cm
\noindent
{\bf Proof:} Suppose that $x=re_0+a$ in $\erre e_0\oplus {\rm Im}(\a_n)$ then $\overline{x}=re_0-a$ and
\begin{eqnarray*}
{\rm exp}(\overline{x})&=&{\rm exp}(re_0-a)=e^r({\rm exp}(-a))=e^r({\rm cos}(||a||)e_0+{\rm sin}(||a||)(-\frac{a}{||a||})\\
&=&e^r(({\rm cos}(||a||)e_0-{\rm sin}(||a||)\frac{a}{||a||})\\
&=&\overline{{\rm exp}(x)}.
\end{eqnarray*}
so we are done with (i).

To prove (ii) recall for non-zero $y$ in $\a_n,$ by definiton, $y^{-1}=||y||^{-2}\overline{y}$  so 
\begin{eqnarray*}
{\rm exp}(x)^{-1}&=&||{\rm exp}(x)||^{-2r}{\rm exp}(\overline{x})\quad {\rm by} (1)\\
&=&e^{-2r}(e^r({\rm exp}(-a)))=e^{-r}{\rm exp}(-a)\\
&=&{\rm exp}(-x)
\end{eqnarray*}
{\bf Corollary 1.9} For $x$ in $\a_n$ and $k$ in $\ze$.
$${\rm exp}(kx)=({\rm exp}(x))^k\quad\hbox{\rm (De Moivre's Formula)}.$$
{\bf Proof:} For $k=0.$ The assertion is trivial.

For $k>0$. The proof is straight-forward using the addition of angles identities, for sine and cosine respectively and induction on $k$.

For $k<0$. The proof follows from case $k>0$ and Proposition 1.8 (ii).

\hfill Q.E.D.

\vskip1cm
\noindent
{\bf II. The exponential map the and $k$-power map.}
\vskip.3cm
$S(\a_n)=S^{2^n-1}$ denotes the unit sphere in $\a_n=\erre^{2^n}$.
\vskip.3cm
\noindent
{\bf Theorem 2.1} The exponential map  exp:$\a_n\rightarrow\a_n\backslash\{0\}$ and its restriction exp:${\rm Im} (\a_n)\rightarrow S(\a_n)$ are onto maps for all $n\geq 1$.

\vskip.3cm
\noindent
{\bf Proof:} By Corollary 1.2 we know that exp (Im $(\a_n))\subset S(\a_n)$. 

Suppose that $y=se_0+b$ for $b$ in Im$(\a_n)$ and $s$ in $\erre$ with $||y||^2=s^2+||b||^2=1$.

If $b=0$ then $s^2=1$ and $s=\mp 1$ and $y=\mp e_0$ but ${\rm exp}(0)=e_0$ and exp$(\pi c)=-e_0$ for all $c\in S(\a_n)$.

Suppose that $b\neq 0$ in Im $(\a_n)$ then there is a real number $\theta$ such that $0<\theta <\pi$ with $s={\rm cos}(\theta)$ and $||b||={\rm sin}(\theta)$.

Let us define $a$ as the non-zero element in Im $(\a_n)$ of norm $\theta$ and linearly dependent to $b$, this means,
$$||a||=\theta\quad{\rm and}\quad||b||(\frac{a}{||a||})=b.$$
Therefore 
\begin{eqnarray*}
{\rm exp}(a)&=&{\rm cos}(||a||)e_0+{\rm sin}(||a||)\frac{a}{||a||}\\
&=&{\rm cos}(\theta)e_0+{\rm sin}(\theta)\frac{b}{||b||}\\
&=&se_0+b\\
&=&y.
\end{eqnarray*}
Therefore exp:Im$(\a_n)\rightarrow S(\a_n)$ is onto.

Now suppose that $y\neq 0$ in $\a_n$ then $||y||^{-1}y$ is in $S(\a_n)$ so there is $a$ in Im$(\a_n)$ such that exp$(a)=||y||^{-1}y$ then
$$||y||{\rm exp}(a)=y.$$
But the real exponent map $e:\erre\rightarrow\erre^+$ is onto, so there is $r$ in $\erre^+$ such that $||y||=e^r$ and if $x=re_0+a$ we have that
$${\rm exp}(x)=e^r{\rm exp}(a)=y$$

\hfill Q.E.D. 

\vskip.5cm
Now we study the $k$-power map $x\mapsto x^k$ for $k$ in $\ze$.
\vglue.5cm
\noindent
{\bf Lemma 2.2} For $k>0$ and $x$ in $\a_n$ we have

(1) $\overline{(x^k)}=(\overline{x})^k$

(2) $||x||^k=||x^k||$

(3) If $x\neq 0$ then $(x^{-1})^k=(x^k)^{-1}$
\vglue.5cm
\noindent
{\bf Proof:} Notice that for $k=1$ (1), (2) and (3) are trivial. So we assume that $k\geq 2$.

(1) We proceed by induction on $k$.

Recall that for $x$ and $y$ in $\a_n$,   $\overline{xy}=\overline{y}\overline{x}$ so
$(\overline{x})^2=(\overline{x})(\overline{x})=\overline{xx}=
\overline{(x^2)}$. Suppose now that $(\overline{x})^k=(\overline{x^k})$. Then $(\overline{x})^{k+1}=(\overline{x})^k(\overline{x})=\overline{(x^k)}\overline{x}=\overline{xx^k}=\overline{x^{k+1}}$ and (1) is done.

(2) Notice that for $x$ and $y$ in $\a_n$ if $\overline{x}(xy)=(\overline{x}x)y=||x||^2y$ then $$||xy||=||x||||y||$$  because
$$||xy||^2=\langle xy,xy\rangle=\langle y,\overline{x}(xy)\rangle=\langle y,||x||^2y\rangle=||x||^2\langle y,y\rangle=||x||^2||y||^2$$
We also notice that $\overline{x}(xy)=(\overline{x}x)y$ if and only if $x(xy)=x^2y$, because $(\overline{x}+x)$ is real 
and hence associates with any other two elements in $\a_n$ so, in particular
$$\overline{x}(xy)+x(xy)=(\overline{x}+x)(xy)=((\overline{x}+x)x)y=(\overline{x}x)y+x^2y=||x||^2y+x^2y$$
therefore $\overline{x}(xy)-||x||^2y=x^2y-x(xy)$.

Making $y=x^{k-1}$ for $k\geq 2$ and recalling that $\a_n$ is power associative, we have that
$||x^{k+1}||=||x(xx^{k-1})||=||(x^2)||||(x^{k-1})||$ for $k\geq 2$. By an obvious induction we are done with (2).

\begin{itemize}
\item[(3)] Recall that for $y\neq 0$ in $\a_n\quad  y^{-1}=||y||^{-2}\overline{y}$.Using (1) and (2) we have that for $x\neq 0$
$$
(x^{-1})^k=(||x||^{-2}\overline{x})^k=||x||^{-2k}(\overline{x})^k =(||x||^k)^{-2}\overline{x}^k)=||x^k||^{-2}(\overline{x}^k)=(x^k)^{-1}
$$
\hfill Q.E.D.
\end{itemize}
\vglue .5cm
\noindent
{\bf Definition:} For nonzero $x$, $x^{-k}:=(x^{-1})^k =(x^k)^{-1}.$

For $k>0\quad \;\rho_k:=\a_n \setminus \{0\}\rightarrow \a_n$ is $\rho_k (x)=x^k.$\\
For $k=0\quad \;\rho_0 :\a_n \setminus \{0\}\rightarrow \a_n$ is $\rho_0 (x)=e_0.$\\
For $k<0\quad \;\rho_k :\a_n \setminus \{0\} \rightarrow \a_n$ is $\rho_k (x)=\rho_{-k}(x^{-1}).$

$\rho_k$ is, by definition, \underline{the $k$--power map} for $k$ in $\ze$; and clearly $\rho_k \circ \rho_\ell=\rho_{k+\ell}$ for $k$ and $\ell$ in $\ze$.
\vglue .5cm
\noindent
{\bf Theorem 2.3} 
For $k\in \ze$ and $\rho_k:S(\a_n)\rightarrow S(\a_n)$
 has topological degree $k$.
\vglue .5cm
\noindent
{\bf Proof:} By Lemma 2.2 (2) we have that    
$\rho_k (S(\a_n))\subset S(\a_n)$. 

Now define $\sigma_k:\a_n \rightarrow \a_n$ as $\sigma_k (x)=kx$ 
for $k\in \ze.$
Now by Corollary 1.9 and Lemma 2.2 (2) and (3)  
exp$(\sigma_k (x))=\rho_k$(exp$(x)$).

Therefore Im $(\a_n)$ can be seen, as the tangent space of $S(\a_n)$ at $x=e_0,$ so degree of $\rho_k$ is $k$.

\hfill Q.E.D.
\vglue.5cm

Now we study the $k$--power map $\rho_k:S(\a_{n+1})\rightarrow S(\a_{n+1})=S^{2^{n+1}-1}$ restrcited to a subsphere of dimension $2^n$ for $n\geq 1$.

Consider the vector subspace of $\a_{n+1}=\a_n \times \a_n$ where the second coordinate is real in $\a_n$ that is 
$$
\a_n \times \a_0:=\{(x,y)r\in \a_n \times \a_n|y=re_0\;\;\hbox{\rm for}\;\;\in \erre\}.
$$
Clearly $\a_n \times \a_0$ is a vector subspace of $\a_{n+1},$ which is closed under conjugation and inverses; i.e., if $(x, re_0)\in \a_n \times \a_0$ then $\overline{(x, re_0)}=(\overline{x},- re_o)\in \a_n \times \a_0$ and for $(x, re_0)\neq (0,0)$ then $(x, re_0)^{-1}\in \a_n \times \a_0$. 
\vglue .5cm
\noindent
{\bf Lemma 2.4} The vector subspace $\a_n \times \a_0$ of $\a_{n+1}$ is closed under $k$-powers for $k$ in $\ze$.
\vglue .25cm
\noindent
{\bf Proof:} Clearly the cases $k=0$ and $k=1$ are obvious. Based on  the above observation,the case $k<0$ follows from the case $k>0$. 

First we check that $\a_n \times \a_0$ is closed under squaring operation 
$$(x, re_0)^2=(x, re_0)(x, re_0)=(x^2-r^2e_0, rx+r\overline{x})=
(x^2-r^2e_0, r (x +\overline{x}))$$ 
but $(x+\overline{x})$ is real so $(x, re_0)^2 \in \a_n \times \a_0$. 

Now we proced by induction on $k$ for $k\geq 2$. 

Suppose that $\alpha \in (\a_n \times \a_0)$ and that $\alpha^k \in (\a_n \times \a_0)$. 

We want to prove that $\alpha^{k+1} \in (\a_n \times \a_0)$. 

Since $(\a_n \times \a_0)$ is vector subspace then $(\alpha^k +\alpha)\in (\a_n \times \a_0)$ and because $(\a_n \times \a_0)$ is closed under the squaring operation we have that $(\alpha^k +\alpha)^2, (\alpha^k)^2=\alpha^{2k}$ and $\alpha^2$ are in $\a_n \times \a_0$. 

But, we can associate powers in $\a_{n+1}$ so
$$
(\alpha^k +\alpha)^2 =\alpha^{2k}+2\alpha^{k+1}+\alpha^2.
$$
Therefore $\alpha^{k+1}\in \a_n \times \a_0$

\hfill Q.E.D.
\vglue .5cm
\noindent

{\bf Remark:} We can prove a more general assertion than the one in lemma 2.4.

Let $V$ be a vector subspace of $\a_n$ and define the following vector subspace of $\a_{n+1}.$ 

$$\a_n \times V=\{(x,v)\in \a_n \times \a_n | v\in V\}$$

Clearly these vector subspace is closed under conjugation and inverses and
$$(x,v)^2 =(x^2-||v||^2 e_0,v(\overline{x}+x))$$

But $(\overline{x}+x)$ is real so $\a_n \times V$ is closed under squares 
by similar argument as
in Lemma 2.4. we have that $\a_n \times V$ is closed under $k$ powers for $k$ in $\ze$.
\vglue.5cm
\noindent

{\bf Corollary 2.5} The $k$--power map $\rho_k :S(\a_{n+1})\rightarrow S(\a_{n+1})$ restricted to $S(\a_n \times \a_0)=S^{2^n}$ has also topological degree $k$ for $k\in \ze$. 
\vglue .25cm
\noindent
{\bf Proof:} Recall that $S(\a_{n+1})=\{(x,y)\in \a_n \times \a_n|||x||^2 +||y||^2 =1\}$ thus $S(\a_n \times \a_0)=\{(x, re_0)\in \a_n \times \a_0|||x||^2 +r^2 =1\}$.

By Lemma 2.4 and Lemma 2.2 (2) $\rho_k (S(\a_n \times \a_0))\subset S(\a_n \times \a_0)$ and by Theorem 2.3 $\rho_k :S(\a_n \times \a_0)\rightarrow S(\a_n \times \a_0)$ has degree $k$. 

\hfill Q.E.D.

\vglue .5cm
\noindent
{\bf Remark:} Based in the previous remark we may prove that the $k$ power map on $S(\a_{n+1})$ restricted to $S(\a_n \times V)$ has also degree $k$.

Therefore any continuos map from $S^m$ to itself is homotopic to a $k$ power map
for all $m\geq 1$,because the homotopy class of any continuos map is determined by
degree (Hopf theorem) and if m is between $2^n$ and $2^{n+1}$ we may choose a vector
subspace $V$ of $\a_n$  such that  $S^m =S(\a_n \times V)$ and the restriction
of $\rho_k :S(\a_{n+1})\rightarrow S(\a_{n+1})$ to $S(\a_n \times V)$ is homotopic to the original map. 
\vglue.5cm
\noindent
{\bf III. Fundamental theorem of algebra for $\a_n\;n\geq 1$.}
\vglue .25cm
\noindent
{\bf Definition:} For $x$ non--zero in $\a_n, x=re_0 +a$ in $\erre e_0 \oplus {\rm Im}(\a_n)$ and $k>0$ \underline{$a\;k$-root of $x$} is
$$
\sqrt[k]{x}=||x||^{1/k}\;\hbox{\rm exp}(\frac{1}{k}a)=e^{r/k}\;\hbox{\rm exp}(\frac{1}{k}a).
$$
By Theorem 2.1 every non--zero element in $\a_n$ has (at least) one $k$--root for $k>0$.

\vglue .5cm
\noindent
{\bf Example:} $x^2+e_0=0$ has infinitely many solutions: every element in in the unit
sphere, $S$(Im$(\a_n)$) is a solution.  
\vglue .5cm
\noindent
{\bf Example:} If $a$ is a non--zero pure element in $\a_n$ then $x^k -a=0$ has exactly $k$ solutions, namely, the set of $k$--roots of $a$. 
\vglue.3cm
Now we want to extend the fundamental theorem of algebra for $\a_1=\ce$ to $\a_n$ for $n>1$.

One direct generalization, can be done, on the polynomials which depend only on one imaginary unit in $\a_n$. That is, we look at the polynomials which are $\ce$-dependent to a given $a$ in
 Im $(\a_n)$ with $||a||=1$.
 So $a$ plays the same role, for this polynomials, as $e_1=i$ plays  for complex polynomials and we have a Fundamental Theorem of Algebra in this situation. 
\vglue .5cm
\noindent
{\bf Lemma 3.1} If $x,y$ and $z$ are non-zero $\ce$-dependent elements in $\a_n$ then
\begin{itemize}
\item[(i)] $xy$ is $\ce$- dependent with $z$.
\item[(ii)] $x^k$ and $y^\ell$ are $\ce$-dependent for $k>0$ and $\ell >0$
\item[(iii)] $(xy)z=x(yz)$.
\end{itemize}
\vglue .25cm
\noindent
{\bf Proof:} If one of the three elements is real then the results (i), (ii) and (iii) are obvious.

Suppose that, the three elements $x,y$ and $z$ are non-real, that is, they have non--zero imaginary part. Also, is easy to see, that on the subset of $\a_n$ consisting of non-real elements, $\ce$--dependence define an equivalence relation.

Write $x=re_0+ta\;y=se_0 +qa$ and $z=ue_0 +va$ where $r,t,s,q,u$ and $v$ are real numbers and $a\in {\rm {Im}} (\a_n)$ with $||a||=1$. 

Now
\begin{itemize}
\item[i)] $xy=(re_0+ta)(se_0+qa)=(rs-tq)e_0 +(rq+ts)a$ and $(xy)$ is \ce-dependent with $z$.  

To show (ii) we notice that $x$ is \ce-dependent with $y^2$,because
 $$y^2=(s^2-q^2)e_0+(s+q)a$$. 

Next we proced by induction on $k$. 

Suppose that $x$ is \ce -dependent with $x^k$ and we want to prove that $x$ is \ce -dependent 
with $x^{k+1}$. 

Now $(x^k+x)$, $(x^k+x)^2$, $(x^k)^2$ and $x^2$ are $\ce$-dependent with $x$ so $(x^k+x)^2-(x^k)^2-x^2=2x^{k+1}$ and $x$ is \ce -dependent with $x^{k+1}$.

Since \ce -dependence is an equivalence relation for non-real elements we are done with (ii).

To prove (iii) recall that (see [5])  the associator $(x,y,z):=(xy)z-x(yz)$ is a tri-linear map that vanish if one of the entries is real so by flexibilty

$(x, y, z)=(re_0+ta, se_0 +qa, ue_0 +va)=tqv(a,a,a)=0$.

\end{itemize}
\hfill Q.E.D.
\vglue .5cm
\noindent
{\bf Definition:} A \underline{complex polynomial} in $\a_n$ of degree $k$  is a continuous function of the form
$$
p(x)=\xi_0 +\xi_1 x +\xi_2 x^2+\cdots +\xi_{k-1}x^{k-1}+x^k
$$
where the coefficents $\xi_i$ are $\ce$--dependent among them and with $x,$  and   
 
$\;i=0,1,\ldots , k-1$.

Notice that every polynomial with real coeficients is a complex polynomial.

\vglue .5cm
\noindent
{\bf Theorem 3.2} Every complex polynomial has at least one root in $\a_n$. 
\vglue .25cm
\noindent
{\bf Proof:} This follows from the Fundamental Theorem of Algebra for $\ce$.

Suppose that $x=re_0 +sa$ where $0\neq s$ and $r$ in $\erre$ and $a\in {\rm Im}(\a_n)$
with $||a||=1$ so $p(x)$ and all the summands in $p(x)$ are in the complex subspace of $\a_n$ generated  by $\{e_0, a\}$, because, Lemma 3.1 (i), (ii) and (iii).

Suppose that $x=re_0,$ that is, $s=0,$ so by definition the polynomial is a real polynomial and it has at least one complex root and we may choose any $a\in S(\a_n)$ to inmerse the polynomial into Span$\{e_0, a\}$. 

\hfill Q.E.D.

Now we use what we know about the topology of the $k$--power map to extend the Fundamental theorem of algebra to a more general type of continuous functions than the complex polynomial
in $\a_n$ 

Before that, we show, that some polynomials  have no roots in $\a_n$. 
\vglue .5cm
\noindent
{\bf Exmaple:} For $n\geq 2$ and non--zero $a$ in Im$(\a_n)$ 
$$
p(x)=ax-xa+e_0
$$
has no roots in $\a_n$.Because every commutator of this form,$[a,x]=ax-xa$ has real part equal to zero.
\vglue .5cm
\noindent
{\bf Definition:} A \underline{generalized polynomial} of degree $k$ on $\a_n$ with $k>0$  is a continuous function $p:\a_n \setminus \{0\}\rightarrow \a_n$ of the form
$$
p(x)=x^k(e_0+g(x)).
$$
where $g(x)$ is a nonconstant continuous function,defined for non-zero elements in $\a_n$, such that $||g(x)||\rightarrow 0$ when $||x||\rightarrow \infty$. 
\vglue .5cm
\noindent
{\bf Proposition 3.3} Every complex polynomial in $\a_n$ for $n\geq 1$ is a generalized polynomial in $\a_n$.
\vglue .25cm
\noindent
{\bf Proof:} Suppose that $\xi$ and $x\neq 0$ are $\ce$--dependet in $\a_n$ then by lemma 3.1 and lemma 1.5 (1) we have that
 $$x^{-k}(\xi x^\ell)=(x^{-k}\xi) x^\ell=(\xi x^{-k})x^\ell =\xi (x^{-k}x^\ell)=\xi x^{-k+\ell}$$ 
for $k$ in $\ze$ and $\ell \geq 0.$

Therefore if $p(x)=\xi_0 +\xi_1 x+\cdots +\xi_{k-1}x^{k-1}+x^k$ is a complex polynomial and $g(x):=x^{-k}(\xi_0 +\xi_1x+\cdots +\xi_k x^{k-1})= \xi_0 x^{-k} +\xi_1 x^{-k+1}+\cdots + \xi_k x^{-1}$ then $||g(x)||\rightarrow 0$ when $||x||\rightarrow \infty$ and $p(x)=x^k (e_0 +g(x))$ by Lemma 3.1 (iii).

\hfill Q.E.D.

\vglue .5cm
\noindent
{\bf Theorem 3.4} (Fundamental theorem of algebra for $\a_n$). Every generalized polynomial has at least one root in $\a_n$,for $n\geq 1$.
\vglue .25cm
\noindent
{\bf Proof:} Given $p(x)$ a generalized polynomial in $\a_n$ define $$\hat{p} :\a_n \cup \{\infty\}=S^{2^n}\rightarrow \a_n \cup \{\infty\}=S^{2^n}$$ with
$$
\hat{p}(x)=\left\{
\begin{array}{lcl}
p(x) & \hbox{\rm if} & ||x||<\infty\\
\infty & \hbox{\rm if} & x=\infty
\end{array}\right.
$$
where $\a_n \cup \{\infty\}$ denotes the one-point compactification of $\a_n =\erre^{2^n}$. 

Making the identification
$$
\a_n \cup \{\infty\}=S^{2^n} =S(\a_n \times \a_0)
$$
where the \underline{line at infinity} is $\{(0, re_0)\in \a_n \times \a_0\}$ so $\hat{p}$ is a continuous map from $S(\a_n \times \a_0)$ to $S(\a_n \times \a_0)$.

{\bf Claim:} $\hat{p}$ and $\rho_k$ the $k$--power map, are homotopic.

Let us define 
\begin{eqnarray*}
F_t (x) &=& x^k (e_0 +(1-t)g(x))\\
F_t (\infty) &=& \infty
\end{eqnarray*}
for $0\leq t\leq 1$. Obviously $F_t$ is continuous on $x$ and $t$ and 
$$F_0 (x)=x^k (e_0 +g(x))=p(x)$$ 
and $F_0 (\infty )=\hat{p}(\infty )=\infty$ and $F_1 (x)=x^k$ and $F_1 (\infty )=\infty$.

Thus $\hat{p}$ and $\rho_k$ are homotopic. 

By corollary 2.5, $\rho_k$ has degree $k$ then $\hat{p}$ has degree $k$ and $\hat{p}$  and $p$ are onto, so for $0\in \a_n$ there is $\alpha $ in $\a_n$ such that $p(\alpha )=0$. 

\hfill Q.E.D.

\newpage
\begin{center}
{\bf References}
\end{center}

\begin{enumerate}
\item[1] L.E.~Dickson. On quaternions and their generalization and the history of the eight square theorem. Annals of Mathematics 20, 155-171 and 297, 1919.
\item[2] S.~Eilenberg-I. Niven. The ``Fundamental theorem of algebra'' for Quaternions. Bulletin of the American Mathematical Society 50 246-248 (1949).
\item[3] Eakin-Sathaye.On automorphisms and derivations of Cayley-Dickson algebras.
Journal of Pure and Applied Algebra 129,263-280 1990.
\item[4] S.H.Khalil-P.Yiu.The Cayley-Dickson algebras.A theorem of Hurwitz and quaternions.Boletin de la Sociedad de Lodz. Vol. XLVIII 117-169 1997.
\item[5] G. Moreno. The zero divisors of the Cayley--Dickson algebras over the real numbers. Boletin de la Sociedad Matem\'atica Mexicana (3) Vo. 4 13-27, 1998.
\item[6] G. Moreno. Alternative elements in the Cayley--Dickson algebras. Preprint CINVESTAV, Mexico 2000.and also available at hopf.math.purdue.edu-pub-Moreno-2001.
\item[7] R. Remmert. Numbers. Graduate Text in Mathematics, 123 Springer--Verlag.
\item[8] R.D. Schafer. On the algebras formed by the Cayley-Dickson process. American Journal of Math. 76, 1954 435-446.
\end{enumerate}

\end{document}